\documentclass{amsart}
\usepackage{amsmath, amssymb,verbatim,appendix}
\usepackage[mathscr]{eucal}
\usepackage{amscd}
\usepackage{amsthm}
\usepackage{stmaryrd}
\usepackage{enumerate}
\usepackage{comment}
\usepackage[latin1]{inputenc} 
\usepackage{tikz}
\usetikzlibrary{shapes,arrows}
\usepackage{cite}
\usepackage{url}
\usepackage{tikz-cd}
\usepackage[colorinlistoftodos]{todonotes}

\makeatletter

\newcommand{\dotminus}{\mathbin{\text{\@dotminus}}}

\newtheorem{theorem}{Theorem}[section]
\newtheorem{lemma}[theorem]{Lemma}
\newtheorem{corollary}[theorem]{Corollary}
\newtheorem*{corollary*}{Corollary}
\newtheorem{proposition}[theorem]{Proposition}
\newtheorem{remark}[theorem]{Remark}

\newtheorem{question}[theorem]{Question}
\newtheorem{conjecture}[theorem]{Conjecture}

\newtheorem{fact}[theorem]{Fact}

\newtheorem{innercustomthm}{Theorem}
\newenvironment{customthm}[1]
  {\renewcommand\theinnercustomthm{#1}\innercustomthm}
  {\endinnercustomthm}

\theoremstyle{definition}

\newtheorem{definition}[theorem]{Definition}

\def\A{\mathbb{A}}

\def\P{\mathbb{P}}

\newcommand{\calF}{\mathcal{F}}

\newcommand{\calC}{\mathcal{C}}

\newcommand{\calQ}{\mathcal{Q}}

\newcommand{\calU}{\mathcal{U}}

\newcommand{\PGL}{\operatorname{PSL}}
\newcommand{\SL}{\operatorname{SL}}
\newcommand{\GL}{\operatorname{GL}}

\newcommand{\Span}{\operatorname{span}}
\newcommand{\dcf}{\operatorname{DCF}}
\newcommand{\ccm}{\operatorname{CCM}}
\newcommand{\acf}{\operatorname{ACF}}
\newcommand{\id}{\operatorname{id}}
\newcommand{\Autdef}{\mathrm{Aut}_{\mathrm{def}}}

\def\mr{\operatorname{RM}}
\def\acl{\operatorname{acl}}
\def\dcl{\operatorname{dcl}}
\def\tp{\operatorname{tp}}

\newcommand{\SO}{\operatorname{SO}}
\newcommand{\PSp}{\operatorname{PSp}}
\newcommand{\PSO}{\operatorname{PSO}}
\newcommand{\Gr}{\operatorname{Gr}}
\newcommand{\Sp}{\operatorname{Sp}}

\def\Ind#1#2{#1\setbox0=\hbox{$#1x$}\kern\wd0\hbox to 0pt{\hss$#1\mid$\hss}
\lower.9\ht0\hbox to 0pt{\hss$#1\smile$\hss}\kern\wd0}
\def\ind{\mathop{\mathpalette\Ind{}}}
\def\Notind#1#2{#1\setbox0=\hbox{$#1x$}\kern\wd0\hbox to 0pt{\mathchardef
\nn=12854\hss$#1\nn$\kern1.4\wd0\hss}\hbox to
0pt{\hss$#1\mid$\hss}\lower.9\ht0 \hbox to
0pt{\hss$#1\smile$\hss}\kern\wd0}
\def\nind{\mathop{\mathpalette\Notind{}}}

\newcommand{\deltatype}{\Delta\mbox{-}\operatorname{type}}

\title{Finite-dimensional differential-algebraic permutation groups}
\author{James Freitag}
\address{James Freitag\\
University of Illinois Chicago\\ 
Department of Mathematics, Statistics,
and Computer Science\\ 
851 S. Morgan Street\\
Chicago, IL, 60607-7045\\
USA}
\email{jfreitag@uic.edu}

\author{L\'eo Jimenez}
\address{L\'eo Jimenez\\
The Ohio State University\\
Department of Mathematics\\
231 West 18th Avenue\\
Columbus, OH \ 43210-1174\\
USA}
\email{jimenez.301@osu.edu}

\author{Rahim Moosa}
\address{Rahim Moosa\\
University of Waterloo\\
Department of Pure Mathematics\\
200 University Avenue West\\
Waterloo, Ontario \  N2L 3G1\\
Canada}
\email{rmoosa@uwaterloo.ca}

\begin{document}

\date{\today}

\thanks{R. Moosa was partially supported by an NSERC Discovery Grant and a Waterloo Math Faculty Research Chair.
L. Jimenez was partially supported by the Fields Institute for Research in Mathematical Sciences.
J. Freitag was partially supported by NSF CAREER award 1945251.}

\begin{abstract}
Several structural results about permutation groups of finite rank definable in differentially closed fields of characteristic zero (and other similar theories) are obtained.
In particular, it is shown that every finite rank definably primitive permutation group is definably isomorphic to an algebraic permutation group living in the constants.
Applications include the verification, in differentially closed fields, of the finite Morley rank permutation group conjectures of Borovik-Deloro and Borovik-Cherlin.
Applying the results to binding groups for internality to the constants, it is deduced that if complete types $p$ and $q$ are of rank $m$ and $n$, respectively, and are nonorthogonal, then the $(m+3)$rd Morley power of $p$ is not weakly orthogonal to the $(n+3)$rd Morley power of $q$.
An application to transcendence of generic solutions of pairs of algebraic differential equations is given.
\end{abstract}

\maketitle

\setcounter{tocdepth}{1}
\tableofcontents

\section{Introduction}

\noindent 
Classical Galois theory connects problems on polynomial equations over fields with purely group-theoretic statements. The connection is especially fruitful because the latter is an arena in which some of the most powerful classification results have been obtained.
In this article, we exploit an analogous connection between differential equations and actions of differential algebraic groups, where we establish several classification results and apply them to algebraic differential equations. 

Specifically, this article is concerned with the structure of differential-algebraic permutation groups of finite rank.
Here, by a {\em permutation group} we mean a group~$G$ acting faithfully and transitively on a set $S$.
By {\em differential-algebraic} we mean that the group, the set, and the action, are all described by algebraic differential equations over a differential field $(k,\delta)$ of characteristic zero.
Equivalently, they are definable in the first order theory of differentially closed fields ($\dcf_0$).
(See~\cite[Chapter~II]{marker2017model} for a detailed introduction to the model theory of algebraic differential equations, and~\cite{moosa-dcf} for a quick one.)
Finally, by {\em finite rank}, we mean that for any $a\in G$, the differential field generated by $a$ over $k$ is of finite transcendence degree.

We are motivated to study the structure of such group actions because they arise spontaneously in the analysis of the algebraic relations between solutions of algebraic differential equations.
An instance of such an application will be explained at the end of this Introduction.

The field of constants in a differentially closed field is a pure algebraically closed field, and this allows us to view finite rank differential-algebraic geometry as an expansion of algebraic geometry.
Much of the work we do here involves comparing differential-algebraic permutation groups with {\em algebraic} permutation groups, namely, algebraic groups acting on algebraic varieties in the field of constants.

What follows is a summary of our results.

\smallskip
\subsection{Definably primitive permutation groups}
A definable permutation group $(G,S)$ is {\em definably primitive} if $S$ admits no definable proper nontrivial $G$-invariant equivalence relations.
Numerous questions in permutation group theory can be reduced to the case of primitive group actions. We show that every such action in finite rank differential-algebraic geometry comes from pure algebraic geometry:

\begin{customthm}{A}
\label{A}
Every connected finite rank definably primitive definable permutation group in $\dcf_0$ is definably isomorphic to the constant points of an algebraic permutation group over the constants.
\end{customthm}

This is Theorem~\ref{dppg-dcf} below, and is proved by using the finite Morley rank O'Nan-Scott type theorem of Macpherson and Pillay~\cite{macpherson1995primitive} to reduce to the case of simple differential-algebraic groups, and then using Cassidy's theorem~\cite{cassidy1989classification} in that case.

\smallskip
\subsection{Base size}
By a {\em base} for a permutation group $(G,S)$ we mean a subset of $S$ whose pointwise stabiliser is trivial.
A conjecture of Borovik and Deloro~\cite{borovik2019binding} predicts that in the finite Morley rank setting, where bases are always finite, the size of a minimal base for a definably primitive permutation group grows at most linearly with the Morley rank of $S$.

Again, though this problem is stated purely in terms of group theoretic aspects of definable group actions, it has a motivation from the perspective of algebraic differential equations. Suppose we are given an equation $X$ of order $n$ whose general solution can be written as a rational function of $m$-many constants $c_1, \ldots , c_m$ and $p$-many solutions $a_1, \ldots , a_p$ of $X$.\footnote{Equivalently, $X$ is internal to the constants.} Such an expression is the key to constructing new solutions of $X$ by a \emph{superposition} of existing solutions, and is an extensively studied topic in differential equations \cite{jones1967nonlinear}. It is natural to ask if some general bound on $p$ can be given in terms of the order (and perhaps the degree) of the equations defining $X$. No such general result seems to exist in the present literature, but the problem can readily be seen to be equivalent to bounding the size of a base for the action of the binding group on $X$. In this section, we give just such a bound, under the additional assumption that the action is definably primitive. The following theorem, which proves the conjecture of Borovik and Deloro in $\dcf_0$, appears as Corollary~\ref{bd-dcf} below:

\begin{customthm}{B}
\label{B}
There is a constant $c$ such that if $(G,S)$ is any connected finite rank definably primitive definable permutation group in $\dcf_0$ then there is a base of size less than $c\mr(S)$.
\end{customthm}

While the theorem is stated here for differential-algebraic permutation groups, Theorem~\ref{A} immediately reduces the situation to the consideration of algebraic permutation groups in the constants.
That is, one needs only to prove the theorem for $\acf_0$, the theory of algebraically closed fields in characteristic zero.
This is done by again carrying out an O'Nan-Scott type analysis and reducing to the case when $G$ is a simple algebraic group, and then applying the results of Burness et. al. from~\cite{burness2017base}. 

As binding group actions need not be definably primitive, Theorem~B does not always apply to the motivating problem about algebraic differential equations discussed above.
We leave that for future work:

\begin{question}
Given an algebraic differential equation $X$ that is internal to the constants, can one bound the size of the base for binding group action in terms of some invariants of $X$?
\end{question}

\smallskip
\subsection{Multiple transitivity and the Borovik-Cherlin Conjecture}
For an integer $\mu>1$, a permutation group $(G,S)$ is {\em $\mu$-transitive} if the co-ordinate-wise action on $S^\mu$ is transitive off the diagonals -- that is if $G$ takes any tuple of $\mu$ distinct elements of $S$ to any other tuple of $\mu$ distinct elements.
One instance of the fact that there are few $\mu$-transitive group actions is Knop's~\cite{knop1983mehrfach} classification of all $2$-transitive algebraic group actions: the only possibilities are  $\PGL_{n+1}$ on $\mathbb P^n$, or certain algebraic subgroups of the group of affine transformations on $\mathbb A^n$, for some $n>1$.
An immediate consequence of Theorem~\ref{A}, once you observe that $2$-transitivity implies primitivity, is an analogous classification of finite rank $2$-transitive permutation groups in $\dcf_0$, see Theorem~\ref{knop} below.

A more flexible notion of multiple transitivity is {\em generic} $\mu$-transitivity.
Here we are in the context of a finite Morley rank permutation group $(G,S)$ and we ask that the co-ordinatewise action of $G$ on $S^\mu$ admits an orbit that is generic in the sense that its complement is of strictly smaller Morley rank than $S^\mu$.
The notion was introduced and studied by Borovik and Cherlin in~\cite{BC2008} as an abstraction of generic transitivity for algebraic groups in the sense of Popov~\cite{Popov2007}, with which it agrees if one is working in $\acf_0$.
There are many more examples of generic $\mu$-transitivity, but a conjecture of Borovik and Cherlin predicts that if $G$ acts generically $(n+2)$-transitively, where $n=\mr(S)$, then $(G,S)$  is isomorphic to the natural action of $\mathrm{PSL}_{n+1}(F)$ on~$\mathbb{P}^n(F)$, for some algebraically closed field~$F$.
In $\acf_0$ this was verified by the first and third authors in~\cite{nmdeg}, following a strategy suggested in~\cite{BC2008}.
The conjecture remains largely open otherwise, the only known case for arbitrary theories is when $\mr(S)=2$, dealt with by Alt\i nel and Wiscons in~\cite{altinel2018recognizing}.
Here we establish the conjecture for finite rank definable group actions in $\dcf_0$.

\begin{customthm}{C}
\label{C}
Suppose $(G,S)$ is a connected definable permutation group in $\dcf_0$ with $G$ of finite rank and $n=\mr(S)>0$.
If the action is generically $(n+2)$-transitive then $(G,S)$ is definably isomorphic to the natural action of the constant points of $\mathrm{PSL}_{n+1}$ on~$\mathbb{P}^n$.
\end{customthm}

This appears as Theorem~\ref{bc} below, and is proved by reducing to the definably primitive case and then applying Theorem~\ref{A}.
In Section~\ref{sect:mt} we also articulate a more geometric version of this theorem where dimension takes the place of Morley rank, and then conjecture that this geometric formulation holds without the finite rank assumption (even in the context of several commuting derivations).
See Theorem~\ref{fdgbc} and Conjecture~\ref{gbc} below.

\smallskip
\subsection{Nonorthogonality and applications}
As we have already mentioned, our motivation for the study of permutation groups in $\dcf_0$ comes from the fact that they arise as binding groups for internality.
An understanding of definable permutation groups can thus contribute to the study of the fine structure of finite rank types, which in turn has applications to algebraic differential equations.
This connection was made in~\cite{nmdeg} and exploited further in~\cite{freitag2022any} and~\cite{abred}.
For example, as a quick consequence of Theorem~\ref{C}, we show in Corollary~\ref{cor:int} below, that if a type $p$ of Morley rank~$n$ is nonorthogonal to a definable set $X$ then the Morley power $p^{(n+3)}$ is not weakly orthogonal to~$X$.
See Section~\ref{nonorth} for details, including a review of the geometric-stability notions involved.
In fact, using a more careful analysis, we show, in Theorem~\ref{orthbound} below:\footnote{Actually, the proof of Theorem~\ref{D} uses only the truth of the Borovik-Cherlin Conjecture in $\acf_0$, and so, as such, is not really an application of the previous results.}

\begin{customthm}{D}
\label{D}
Suppose $p$ and $q$ are complete stationary types in $\dcf_0$ of $U$-rank $m$ and $n$, respectively.
If $p$ is nonorthogonal to $q$ then $p^{(m+3)}$ is not weakly orthogonal to $q^{(n+3)}$.
\end{customthm}

These results accomplish, for finite rank types in $\dcf_0$, the general goals  mentioned by Hrushovski in~\cite[Section~2, first paragraph]{hrushovski1989almost}.

We now explain an application of Theorem~\ref{D} to the transcendence of solutions of algebraic differential equations.
Fix a characteristic zero algebraically closed differential field $(k,\delta)$ and consider an order~$n$ differential equation:
$
P(y,\delta y,\dots,\delta^{(n)}y)=0$,
where $P\in k[x_0,\dots,x_n]$ is irreducible.
For each $m\geq 1$, consider the following condition on this equation:
\begin{itemize}
\item[$(C_m)$]
For any $m$ distinct solutions $a_1,\dots,a_m\notin k$ the sequence
$$(\delta^{(i)}a_j:i=0,\dots,n-1,j=1,\dots,m)$$
is algebraically independent over $k$.
\end{itemize}
In \cite{freitag2022any}, it was shown that $(C_3)$ implies $(C_m)$ for all $m$, assuming that $n\geq 2$.
That is, to detect whether there are algebraic relations between solutions of a differential equation and its derivatives up to the order of the equation, it surprisingly suffices to consider only triples of solutions. This result, along with a stronger form proved in \cite{nmdeg<2}, has recently been applied to prove new transcendence results for functions satisfying various differential equations \cite{casale2022strong, devilbiss2021generic, freitag2017algebraic}.

Next, consider the variant of the previous condition involving algebraic relations between generic solutions of \emph{different equations}. 
Consider two algebraic differential equations of order $n_1$ and $n_2$, respectively: 
\begin{eqnarray}
\label{ode2} P_1(y,\delta y,\dots,\delta^{(n_1)}y)=0\\
\label{ode3} P_2(y,\delta y,\dots,\delta^{(n_2)}y)=0
\end{eqnarray}
where $P_i\in k[x_0,\dots,x_{n_i}]$ are irreducible, for $i=1,2$.
For each $m_1, m_2\geq 1$, consider the following condition:
\begin{itemize}
\item[$(C_{m_1,m_2})$]
For any $m_1$ independent generic solutions $a_1,\dots,a_{m_1}$ of~(\ref{ode2}) and $m_2$ independent generic solutions $b_1,\dots,b_{m_2}$ of~(\ref{ode3}), the sequence $(\delta^{(i)}a_j, \delta^{(h)}b_\ell)$ where $i=0,\dots,n_1-1$, $h=0,\dots, n_2-1$, $j=1,\dots m_1$, and $\ell=1,\dots, m_2$,
is algebraically independent over $k$.
\end{itemize}
That is, the condition $(C_{m_1,m_2})$ fails when there are nontrivial algebraic relations between $m_1$ generic solutions of equations (\ref{ode2}) and $m_2$ generic solutions of (\ref{ode3}), along with their derivatives up to the order of the equations.
A consequence of Theorem~\ref{D} is the following:

\begin{corollary*}
\label{cor:ade}
$(C_{n_1+3,n_2+3})$ implies $(C_{m_1,m_2})$ for all $m_1,m_2$.
\end{corollary*}

\begin{proof}
Let $p,q\in S(k)$ be the Kolchin-generic types of~(\ref{ode2}) and~(\ref{ode3}), respectively.
There are unique such generic types by irreducibility of the $P_i$.
The realisations of $p^{(\ell)}$ are precisely the sequences of $\ell$ generic independent solutions to~(\ref{ode2}), and similarly for $q$.
And the failure of $(C_{m_1,m_2})$ is precisely dependence between a realisation of $p^{(m_1)}$ and a realisation of $q^{(m_2)}$, that is, $p^{(m_1)}$ and $q^{(m_2)}$ being not weakly orthogonal.
Let $m=U(p)\leq n_1$ and $n=U(q)\leq n_2$.
Hence, if $(C_{m_1,m_2})$ fails then $p^{(m_1)}$ and $q^{(m_2)}$ are not weakly orthogonal, which in turn implies that $p$ and $q$ are nonorthogonal, which by
Theorem~\ref{D} implies that $p^{(m+3)}$ and $q^{(n+3)}$ are not weakly orthogonal, so that $(C_{m+3,n+3})$ fails, and hence $(C_{n_1+3,n_2+3})$ fails.
\end{proof}

\smallskip
\subsection{Beyond $\dcf_0$}
While we prove things exclusively for $\dcf_0$, our results hold in the case of several commuting derivations ($\dcf_{0,m}$), and also in the theory of compact complex manifolds ($\ccm$).
We point this out in a final Section~\ref{ot}, where we articulate abstractly the conditions on a first order theory that we require for our proofs to go through.

\bigskip
\section{Definably primitive permutation groups in $\dcf_0$}

\noindent
Recall that a definable group action $(G,S)$ is said to be {\em definably primitive} if $S$ admits no definable proper nontrivial $G$-invariant equivalence relations.
When the action is transitive, and so $S=G/H$ for some definable group $H\leq G$ with the action being left multiplication, definable-primitivity is equivalent to $H$ being a maximal proper definable subgroup of $G$.

The purpose of this section is to show that there are no new finite rank definably primitive permutation groups in differentially closed fields; they all come from algebraic groups in the constants. We work in a saturated model $(\mathcal{U},\delta) \models \mathrm{DCF}_0$ that will serve as a universal domain for differential-algebraic geometry.
We denote by~$\calC$ the field of constants.

We will make use of the following structure theorem for simple differential algebraic groups of finite rank.
The theorem is originally due to Cassidy~\cite{cassidy1989classification} without the finite rank assumption, with an easier proof in the finite rank case given by Pillay in~\cite[Chapter~III, Theorem~1.5]{marker2017model}.

\begin{fact}
\label{cassidy-pillay}
Every simple finite rank definable group in $\dcf_0$ is definably isomorphic to the $\calC$-points of a simple linear algebraic group over~$\calC$.
\end{fact}

The following is likely well known.

\begin{lemma}
\label{nosbgp}
Suppose $G$ is a definable group of finite rank in $\dcf_0$ that has no proper nontrivial definable subgroups.
Then $G$ is definably isomorphic to the $\calC$-points of an algebraic group over~$\calC$.
\end{lemma}

\begin{proof}
Note that $G$ is connected, and we may assume it is nontrivial.
When $G$ is noncommutative, then, as it is definably simple, it is outright simple (see \cite[Corollary 5.9]{wagner1997stable}), and the result follows from Fact~\ref{cassidy-pillay} above.

So we may assume that $G$ is commutative.
In this case, we show that $G$ is definably isomorphic to $\mathbb G_a(\calC)$.
We can embed $G$ as a Zariski-dense definable subgroup of a (necessarily commutative) algebraic group $E$, see~\cite[Chapter~III, Lemma~1.1]{marker2017model}.
Let $E$ be such of minimal dimension.
As $G$ has no proper nontrivial definable subgroups, it has trivial intersection with any proper nontrivial algebraic subgroup of~$E$.
Modding out by any proper infinite algebraic subgroup would therefore embed $G$ into a smaller-dimensional algebraic group, contradicting the choice of minimal dimension.
So $E$ has no proper infinite algebraic subgroups.
The only possibility for such is $\mathbb G_a, \mathbb G_m,$ or a simple abelian variety.
Now, by~\cite[Lemma~2.5]{hrushovskiminimal}, every infinite definable subgroup of a simple abelian variety contains the Manin kernel, and hence in particular has nontrivial finite subgroups coming from torsion, so $G$ cannot embed in a simple abelian variety.
If $E$ is $\mathbb G_m$ then the logarithmic derivative map restricts to a definable homomorphism $\phi:G\to\mathbb G_a$.
The kernel of $\phi$ must either be trivial or all of $G$, but the latter is impossible as in that case $G=\mathbb G_m(\calC)$ which again has nontrivial finite subgroups.
So $\phi$ embeds $G$ in $\mathbb G_a$, and we may assume that $E=\mathbb G_a$.
Every definable subgroup of $\mathbb G_a$ is a $\calC$-vector subspace.
That $G$ has no proper nontrivial definable subgroups implies it is $1$-dimensional.
Hence, $G$ is definably isomorphic to $\mathbb G_a(\calC)$, as desired.
\end{proof}

Our study of finite rank definably primitive group actions in $\dcf_0$ will rely heavily on the O'Nan-Scott analysis of primitive permutation groups of finite Morley rank carried out by Macpherson and Pillay~\cite{macpherson1995primitive}.
The following fact summarises the aspects of that analysis that we will make use of explicitly; it is drawn from both the statements and proofs of~\cite[Theorems~1.1 and~1.2]{macpherson1995primitive}.

\begin{fact}
    \label{fact:onanscott}
    Suppose $G$ is a connected group of finite Morley rank, and $H$ is a proper definable subgroup such that the action of $G$ on $G/H$ is faithful and definably primitive.
    Let $B$ be the {\em definable socle} of $G$, namely the subgroup generated by the minimal normal definable subgroups of $G$.
    Then $B$ itself is a normal definable subgroup, and one of the following cases holds:
    \begin{enumerate}
        \item $B$ is elementary abelian.
        \item $B$ is torsion-free divisible abelian and acts regularly on $G/H$.
        Moreover, if $H$ is infinite, then there is a finite rank definable algebraically closed field $K$ such that $B$ has a definable finite dimensional $K$-vector space
structure and the action of $H$ on $B$ by conjugation gives an embedding of~$H$ into $\GL(B,K)$.
\item $B$ is the unique minimal normal definable subgroup of $G$, it has trivial centraliser in $G$, and it is simple.
\item $G$ has exactly two minimal definable normal subgroups $T_1$ and $T_2$, both simple, both acting regularly on $G/H$, with $C_G(T_1) = T_2$ and $C_G(T_2) = T_1$, and such that $B$ is the direct product of $T_1$ and $T_2$.
    \end{enumerate}
\end{fact}

\begin{remark}
Let us clarify how our presentation of the O'Nan-Scott analysis in Fact~\ref{fact:onanscott} above compares to the statements in~\cite{macpherson1995primitive} .
Cases~(1) and~(2) correspond to Case~1 of~\cite[Theorem~1.1]{macpherson1995primitive} taking into account also ~\cite[Theorem~1.2]{macpherson1995primitive}.
Our Case~(3) corresponds to Cases~2 and~3 of~\cite[Theorem~1.1]{macpherson1995primitive}.
As explained in the discussion following~\cite[Theorem~1.1]{macpherson1995primitive},  connectedness rules out Cases~4(a)(i) and~4(b) of that theorem.
Our Case~(4) thus corresponds to Case~4(a)(ii) of~\cite[Theorem~1.1]{macpherson1995primitive}.
Finally, the fact that in Case~(3) the centraliser of~$B$ is trivial, and that in Case~(4) we have $C_G(T_1) = T_2$ and $C_G(T_2) = T_1$, follows from the argument in the first two paragraphs of the proof of Theorem~1.1 on page~487 of~\cite{macpherson1995primitive}.
\end{remark}

\begin{theorem}
\label{dppg-dcf}
Suppose $G$ is a connected definable group of finite rank in $\dcf_0$.
If $G$ acts definably, faithfully, transitively, and definably primitively on some infinite set, then $G$ is definably isomorphic to the constant points of an algebraic group over the constants.
\end{theorem}

\begin{proof} 
Taking $H$ to be the point stabilizer of some point in the set on which $G$ is acting, we see that we are in the situation of Fact~\ref{fact:onanscott}.
That is, $G$ acts faithfully and definably primitively on $G/H$.
In particular, $H$ is a maximal proper definable subgroup of $G$.
Let $B$ be the definable socle of $G$.
We proceed by analysing the different possibilities for~$B$ enumerated in Fact~\ref{fact:onanscott}.

As no infinite elementary abelian groups are definable in $\dcf_0$, Case~(1) of Fact~\ref{fact:onanscott} cannot occur.

Suppose we are in Case~(2).
In particular, $B$ acts regularly on $G/H$.
This means that $G=BH$ and $B\cap H=(1)$.
That is, $G=B\rtimes H$ is the semidirect product of $B$ by $H$.
Here $H$ acts naturally on $B$ by conjugation.

If $H$ is finite then the connectedness of $G$, and the existence of a definable surjective homomorphism $\pi:G\to H$ with $B=\ker(\pi)$, imply that $H=(1)$.
As $H$ was maximal, this means that $G$ has no nontrivial proper definable subgroups.
Lemma~\ref{nosbgp} then tells us that $G$ is definably isomorphic to the constant points of an algebraic group over the constants, as desired.

So we may assume that $H$ is infinite, and so the "moreover" clause of Case~(2) applies.
That is, there is a finite rank definable algebraically closed field $K$ such that $B$ has a definable finite dimensional $K$-vector space
structure and the action of $H$ on $B$ by conjugation gives an embedding of~$H$ into $\GL(B,K)$.
As the only finite rank infinite definable field in $\dcf_0$, up to definable isomorphism, is the field of constants, we may assume that $K=\calC$.
Fixing a basis for $B$ over $\calC$, we obtain a definable isomorphism $\nu:B\to\calC^n$ and a definable embedding $\rho:H\to\GL_n(\calC)$, inducing a definable embedding of $G=B\rtimes H$ into $\calC^n\rtimes\GL_n(\calC)$.
By stable embeddedness of the constants, the image of G will then itself be the constant points of an algebraic group.
This completes the proof in Case~(2).

Suppose we are in Case~(3), where $B$ is the unique minimal normal definable subgroup of $G$, has trivial centraliser in $G$, and is simple.
By Fact~\ref{cassidy-pillay}, simplicity implies that $B$ is definably isomorphic to $T:=E(\calC)$ where $E$ is a simple linear algebraic group over $\calC$.
It will therefore suffice to show that in this case $G=B$.
As $G$ is connected, it suffices to show that $G/B$ is finite.
It is easily verified that any definable isomorphism $\phi:B\to T$ extends to an (abstract) isomorphism
    $$\hat\phi:\Autdef(B)\to\Autdef(T)$$
given by $\hat\phi(f)=\phi f\phi^{-1}$.
Here, as $B$ is simple, we can (and do) identify $B$ with its group of inner automorphisms in $\Autdef(B)$, and similarly for $T$.
It follows that $\Autdef(B)/B$ is isomorphic to $\Autdef(T)/T$.
Now, the definable automorphisms of $T$ are just the algebraic automorphisms of $T$ viewed as an algebraic group in the constants.
As $T$ is a simple linear algebraic group, $\Autdef(T)/T$ is finite, see~\cite[$\S$27.4]{humphreys}.
So it remains to observe that $G$ embeds into $\Autdef(B)$ over $B$ via the action by conjugation -- which follows from the fact that $B$ has trivial centraliser in $G$.

Finally, suppose we are in Case~(4) of Fact~\ref{fact:onanscott}.
That is $G$ has exactly two minimal definable normal subgroups $T_1$ and $T_2$, both simple, both acting regularly on $G/H$, and $B$ is the direct product of $T_1$ and $T_2$. Moreover $C_G(T_1) = T_2$ and $C_G(T_2) = T_1$, thus $C_G(B)$ is trivial. We will again show that $G=B$, which will suffice as each $T_i$ is definably isomorphic to the constant points of a simple linear algebraic group in the constants. Again by connectedness of $G$, it is enough to show that $G/B$ is finite.

Consider the homomorphism $G \to \Autdef(T_1)\times\Autdef(T_2)$ given by:
$$g\mapsto([g]_{T_1},[g]_{T_2})$$
where $[g]_{T_i}$ and $[g]_{T_2}$ is conjugation by $g$ on $T_i$.
Let
$$\rho : G \to \Autdef(T_1)/T_1\times\Autdef(T_2)/T_2$$
be the composition with the quotient map.
Given $g \in G$, we have that
\begin{align*}
    g \in \ker(\rho) & \iff \text{ there is } (t_1,t_2) \in B \text{ with } [g]_{T_1} = [t_1]_{T_1} \text{ and } [g]_{T_2} = [t_2]_{T_2} \\
    & \iff (t_1,t_2)g^{-1} \in C_G(B) \text{ for some } (t_1,t_2) \in B\\
    & \iff (t_1,t_2)g^{-1} = \id_G \text{ for some } (t_1,t_2) \in B\\
    & \iff g \in B \text{ .}
\end{align*}
We thus have a definable embedding of $G/B$ into $\Autdef(T_1)/T_1\times\Autdef(T_2)/T_2$. The latter is finite as $T_1$ and $T_2$ are definably isomorphic to simple linear algebraic groups in the constants. Hence $G/B$ is finite, as desired.
\end{proof}

\bigskip
\section{Base size for primitive permutation groups}

\noindent
Since the constants form a pure algebraically closed field, an expected use of Theorem~\ref{dppg-dcf} would be to extend results about definably primitive permutation groups from $\acf_0$ to $\dcf_0$. In this section we do so for results about the size of bases.

\begin{definition}
A {\em base} for a faithful group action $(G,S)$ is a subset $B\subseteq S$ such that, for any $g\in G$, if $g\upharpoonright_B=\id_B$ then $g=1$.
The minimal cardinality for a base is denoted by $b(G,S)$.
\end{definition}

In~\cite[$\S$2]{borovik2019binding}, Borovik and Deloro predict that, for definably primitive permutation groups of finite Morley rank, there is a linear relationship between base size and rank:

\begin{conjecture}[Borovik-Deloro]
\label{conj:bd}
There is an absolute constant $c$ such that
$$b(G,S) < c\mathrm{RM}(S)$$
for $(G,S)$ any connected definably primitive permutation group of finite Morley rank, with $S$ infinite.
\end{conjecture}

It does not seem that the conjecture has even been verified for algebraic group actions, and that is our aim here. 
Once that is established, Theorem~\ref{dppg-dcf} will imply that the conjecture is true of all finite rank differential-algebraic groups as well.

We will build on the fact that the base size for {\em simple} algebraic group actions has been thoroughly investigated in~\cite{burness2017base}.
First, a remark on terminology: in the literature, what is often meant by a ``simple algebraic group" is an algebraic group that does not contain any proper closed {\bf connected} normal subgroup.
For clarity we will refer to such algebraic groups as \emph{almost simple}, and reserve the term {\em simple algebraic group} for algebraic groups with no proper nontrivial normal algebraic subgroups at all.
Almost simple algebraic groups have finite center, and when we mod out by the center they become simple algebraic groups (even simple as abstract groups).
There is a well known classification theorem stating that every simple algebraic group is either a {\em classical group}, namely belonging to one of the following infinite families: $\mathrm{PSL}_n, \mathrm{PSp}_{2n} , \mathrm{PSO}_{2n}$ and $\mathrm{SO}_{2n+1}$; or is one of finitely many exceptional groups.
The first three classical groups are quotients of the almost simple groups $\mathrm{SL}_{n}, \mathrm{Sp}_{2n},\mathrm{SO}_{2n}$, respectively, by their finite centers.
Note also that because we work over an algebraically closed field $\mathrm{PSL}_n$ is the same as $\mathrm{PGL}_n$.

A careful inspection of the main theorems of Burness et. al. in~\cite{burness2017base} yields the Borovik-Deloro conjecture for simple algebraic group actions.
This is mentioned in~\cite{borovik2019binding}, but we give some explanations:

\begin{proposition}\label{pro: linear-bound-for-simple-alg}
There is a constant~$c$ such that if $G$ is a simple algebraic group acting algebraically on a positive-dimensional variety $S$, in characteristic zero, and such that $(G,S)$ is a definably primitive permutation group, then $b(G,S) < c\dim(S)$.
\end{proposition}

\begin{proof}
We show how this follows from the results in~\cite{burness2017base}.
The group $G$ is either a classical group in a \emph{subspace action} (more details on this case shortly) or not.
If not, then \cite[Theorem 1]{burness2017base} gives that $b(G,S) \leq 6$.
Thus we can now assume that $(G,S)$ is a classical group in a subspace action.
This means that we are in one of the following cases (see, for example, Section~4.1 of~\cite{burness2017base}):
\begin{itemize}
    \item[(1)]
    $G=\PGL_n$ and the action is the one induced by the natural transitive action of $\SL_n$ on $\Gr(n,d)$, the set of $d$-dimensional subspaces of $\A^n$ for some $d\leq\frac{n}{2}$.
    Note that the center is in the kernel of that action, so we do get an induced action.
    \item[(2)]
    $G=\PSp_n$ with $n$ even, and the action is the one induced by the natural transitive action of $\Sp_n$ on either $\mathrm{TS}(n,d)$, the set of totally singular $d$-dimensional subspaces of $\A^n$, or on $\mathrm{ND}(n,d)$, the set of non-degenerate $d$-dimensional subspaces of $\A^n$, for some positive $d\leq\frac{n}{2}$.
    Total singularity and non-degeneracy, here, are with respect to a fixed non-degenerate alternating bilinear form, see~\cite[Chapter~11]{roman2005advanced} for the definitions.
    In both cases, the center of $\Sp_n$ is in the kernel of the action.
    \item[(3)]
    $G$ is either $\PSO_n$ with $n$ even, or $\SO_n$ with $n$ odd, and the action is the one induced by the natural transitive action of $\SO_n$ on either $\mathrm{TS}(n,d)$ or on $\mathrm{ND}(n,d)$, for some positive $d\leq\frac{n}{2}$.
    One works here with respect to a fixed non-degenerate symmetric bilinear form.
    Again, the center is in the kernel of the action.
\end{itemize}
We are using here that if $G$ is a symplectic or orthogonal group then primitivity implies, in characteristic~$0$, that the only subspace actions that appear are totally singular or non-degenerate with respect to the relevant underlying form; see the discussion preceding Theorem~4 in~\cite{burness2017base}.
For more details on totally singular and non-degenerate subspaces, we suggest~\cite[Chapter~11]{roman2005advanced}.
In particular, one can deduce from the information there that we do have the above actions and that they are transitive.

Theorem~4 of~\cite{burness2017base} tells us the base size in each of the above cases, often broken up into further subcases, in terms of $\frac{n}{d}$.
A careful inspection of the statement reveals that,
if $n \geq 7$, then $b(G,S)<\frac{n}{d}+6$ in all cases.
For our purposes it is fine to restrict to $n \geq 7$ as doing so only excludes finitely many group actions.

Thus, to prove the result, it is enough to prove that $$\frac{\frac{n}{d} +6 }{\dim S}$$ is bounded by an absolute constant in each of cases~(1) through~(3).
In case~(1), $S$ is the grassmanian $\Gr(n,d)$ which has dimension $d(n-d)$.
Using that $d\leq \frac{n}{2}$, it is easy to compute that $\displaystyle\frac{\frac{n}{d}+6}{d(n-d)}\leq 14$.
The same computation works for cases~(2) and~(3) if $S$ is $\mathrm{ND}(n,d)$, since the latter is a Zariski open subset of $\Gr(n,d)$ and hence is also of dimension $d(n-d)$.
Finally, we deal with cases~(2) and~(3) with $S=\mathrm{TS}(n,d)$.
The dimension of $\mathrm{TS}(n,d)$ is bounded from below by
$$2d\left(\frac{n-1}{2}-d\right) + \frac{d(d-1)}{2}$$
(see, for example, Section 2 of \cite{li2016equivariant}).
Another easy computation then gives a bound of $16$ for $\displaystyle\frac{\frac{n}{d} + 6}{\dim(S)}$.
\end{proof}

To pass from simple algebraic groups to arbitrary algebraic groups we will apply the O'Nan-Scott analysis.
This will require some information about centralisers of subsets of simple algebraic groups, that we now record:

\begin{fact}
\label{cdim-fact}
There is a constant~$e$ such that, for any infinite simple algebraic group,~$G$, the length of the longest strictly descending chain of centralisers of subsets of $G$ is less than $e\dim G$.
\end{fact}

\begin{proof}
The length of the longest strictly descending chain of centralisers of subsets of $G$ is called the {\em $c$-dimension} and is denoted by $\dim_c(G)$.
See~\cite{shumyatsky2004discriminating} for details on $c$-dimension.
In particular, it isn't hard to see that $\dim_c(\GL_n) = n^2+1$ (see the proof of Proposition 2.1 in~\cite{shumyatsky2004discriminating}, for example), and that the $c$-dimension of a subgroup is at most that of the ambient group (see \cite[Lemma~2.2]{shumyatsky2004discriminating}).

We are free to ignore the finite set of exceptional simple algebraic groups.
As the center must be in any centralizer, taking a central extension does not change the $c$-dimension.
So we may assume that $G$ is either $\operatorname{SL}_{n},$ $\operatorname{Sp}_{n}$ with $n$ even, or $\SO_{n}$.
These groups embed in $\GL_{n}$ and hence have $c$-dimension bounded above by $n^2+1$.
On the other hand, the dimension of $G$ is
$n^2-1, \frac{n(n+1)}{2}, \frac{n(n-1)}{2}$, respectively.
We see that $e=6$ works.
\end{proof}

With these ingredients in place, we can establish the Borovik-Deloro conjecture for algebraic group actions:

\begin{theorem}
\label{bd-acf}
    There is a constant $c$ such that if $(G,S)$ is a connected definably primitive definable permutation group in $\acf_0$, and $S$ is infinite, then
    $$b(G,S) < c\mr(S).$$
\end{theorem}

\begin{proof}
We work in a sufficiently saturated $(\calC,0,1,+,-,\times)\models\acf_0$. We will again use Fact \ref{fact:onanscott}. We can still rule out case~(1). Examining the proof of Theorem~\ref{dppg-dcf} above, we see that carrying out the O'Nan-Scott type analysis of Macpherson and Pillay, but just in $\acf_0$ this time, leads to four following possibilities for $(G,\calC)$, up to definable isomorphism. We list them using the indexing from Fact \ref{fact:onanscott}:

    \begin{enumerate}
        \item[(2)]
        \begin{enumerate}[(a)]
            \item (with finite point stabiliser):
        $\mathbb{G}_a(\mathcal{C})$ acting on itself.
        \item Case~(2) (with infinite point stabiliser):
        $S$ is a finite dimensional $\calC$-vector space, $H\leq\GL(B,\calC)$, and $G=S\rtimes H$ with the natural action on $S$ by affine transformations.
        \end{enumerate}
        \item[(3)]
        $G$ is a simple algebraic group,
        \item[(4)]
        $G=T_1 \times T_2$ where $T_1$ and $T_2$ are simple algebraic groups whose induced actions on $S$ are regular.
    \end{enumerate}

    In~(2)(a) we have that $b(G,S)=1$ so that $c=2$ works.
    
    In~(2)(b) we can take as a base for $(G,S)$ any $\calC$-basis for~$S$ along with the zero vector, so that $b(G,S)$ is at most $\dim_{\calC}(S)+1=\mr(S)+1$.
    Hence $c=3$ works.

    Proposition~\ref{pro: linear-bound-for-simple-alg} deals with~(3), noting that $\mr(S)=\dim(S)$.

Finally, consider Case~(4).
Write $S=G/H=(T_1\times T_2)/H$ where $H$ is a maximal proper definable subgroup.

We claim, first of all, that $H$ is the graph of a definable isomorphism.
Indeed, let $\pi_1:H\to T_1$ and $\pi_2 : H \to T_2$ be the projections.
Surjectivity of the $\pi_i$ follow from regularity of the action of $T_i$ on $G/H$.
Indeed, given $s\in T_1$ let $t\in T_2$ be such that $(1,t)$ takes $(s,1)H$ to $H$; this forces $(s,t)\in H$, showing that $\pi_1$ is surjective.
Similarly $\pi_2$ is surjective.
For injectivity of $\pi_i$, observe that the surjectivity of $\pi_2$ implies that $\pi_2(\ker(\pi_1))$ is a normal subgroup of $T_2$, and hence, by simplicity, is either trivial or all of $T_2$.
If $\pi_2(\ker(\pi_1)) = T_2$ then the surjectivity of $\pi_1$ would imply that 
$H = T_1 \times T_2 = G$, a contradiction.
It follows that $\ker(\pi_1)$ is trivial, and hence  $\pi_1:H\to T_1$ is an isomorphism.
Similarly, $\pi_2:H\to T_2$ is an isomorphism. 
So, we have established that $H = \Gamma(\sigma)$, where $\sigma = \pi_1^{-1} \circ \pi_2 :T_1\to T_2$ is a definable isomorphism.

Next, note that the stabiliser of a coset $gH$ in $G$ is precisely the conjugate $H^g$ of $H$ by $g$.
But in our case, where $G=T_1\times T_2$, and $H=\Gamma(\sigma)$, every coset of $H$ is represented by something of the form $(t,1)$ where $t\in T_1$.
Now
$$H^{(t,1)}=\{(r^t,\sigma(r)):r\in T_1\}$$
and hence
$$H\cap H^{(t,1)}=\{(s,\sigma(s)):s\in C(t)\}$$
where $C(t)$ is the centraliser of $t$ in $T_1$.
In particular, if $t_1,\dots,t_\ell\in T_1$ are such that $\bigcap_{i=1}^\ell C(t_i)=(1)$ then $\{H,(t_1,1)H,\dots,(t_\ell,1)H\}$ is a base for $(G,S)$.
Now, note that if $d=\dim_c(T_1)$ then there exists $t_1,\dots,t_{d-1}\in T_1$ such that $\bigcap_{i=1}^{d-1} C(t_i)=(1)$.
It follows that $b(G,S)$ is at most $\dim_c(T_1)$.
By Fact~\ref{cdim-fact}, $\dim_c(T_1)<e\dim(T_1)$ for some absolute constant $e$.
Since $\mr(S)=\dim(T_1)$ in this case, we have that $c=e$ works here.
\end{proof}

\begin{remark}
It is, of course, natural to ask what the absolute constant $c$ given by Theorem~\ref{bd-acf} is.
Our proof shows that, at least for algebraic groups of sufficiently high dimension, $c=16$ works.
For a complete answer to this question we would need to compute, (1) the base size for certain low-dimensional classical simple algebraic groups actions; namely when $n<7$ and $G$ is $\operatorname{PSL}_n, \PSp_n, \PSO_n$, or $\SO_n$,
and (2) the maximal length of chains of centralizers in the exceptional groups.
\end{remark}

Since there are no new finite rank definably primitive group actions in $\dcf_0$ beyond those in $\acf_0$, by Theorem~\ref{dppg-dcf}, we obtain Conjecture~\ref{conj:bd} for free in $\dcf_0$:

\begin{corollary}
\label{bd-dcf}
There is a constant $c$ such that if $(G,S)$ is a finite rank connected definably primitive definable permutation group in $\dcf_0$, and $S$ is infinite, then
    $$b(G,S) < c\mr(S).$$
\end{corollary}

\begin{proof}
Work in a sufficiently saturated model $(\calU,\delta)\models\dcf_0$ with field of constants~$\calC$.
By Theorem~\ref{dppg-dcf}, $G$ is definably isomorphic to the $\calC$-points of an algebraic group over $\calC$.
Since $S=G/H$ for some definable subgroup $H\leq G$, it follows that $(G,S)$ is definably isomorphic to the $\calC$-points of an algebraic group action over $\calC$.
Now apply Theorem~\ref{bd-acf}.
\end{proof}

\medskip
\subsection{An aside on genericity of minimal bases}
Borovik and Deloro also conjectured, in~\cite[$\S$2]{borovik2019binding}, that, for a definably primitive finite Morley rank permutation group $(G,S)$,  the set of minimal bases is generic in $S^{b(G,S)}$.
Put another way, and borrowing notation from Burness et. al.~\cite{burness2017base}, if we let $b^1(G,S)$ be the smallest positive integer $n$ such that $S^n$ contains a generic subset $U$ with every $n$-tuple in $U$ a base for $G$, then Borovik and Deloro are conjecturing that  $b^1(G,S)=b(G,S)$.
However, already in the context of algebraic permutation groups, Burness et. al. show that in many cases $b^1 (G,S)> b (G,S)$, see for example Theorems 4(ii), 4(iii), 5(iii), 7(iii), and 8 of~\cite{burness2017base}.
We suggest the following amendment:

\begin{conjecture}
\label{conj:bd-revised}
There is a natural number $d$ so that whenever $(G,S)$ is a definably primitive permutation group of finite Morley rank, $b^1 (G,S) - b (G,S) <d.$ 
\end{conjecture} 

A close examination of the results of~\cite{burness2017base} reveals that Conjecture~\ref{conj:bd-revised} holds for actions of simple algebraic groups.

\bigskip
\section{Multiply transitive permutation groups}
\label{sect:mt}
\noindent
Recall that, for $\mu>1$, a definable group action $(G,S)$ is {\em $\mu$-transitive} if the co-ordinatewise action of $G$ on $S^\mu$ is transitive off the diagonal.
And in the finite rank setting, $(G,S)$ is {\em generically $\mu$-transitive} if $(G,S^\mu)$ is transitive off a subset of rank strictly smaller than $\mu\mr(S)$.
In this section we investigate consequences of Theorem~\ref{dppg-dcf} to multiply transitive, and especially generically-multiply transitive, permutation groups definable in $\dcf_0$.

Regarding outright multiple transitivity, there is a complete classification of the  $2$-transitive algebraic group actions.
As $2$-transitivity implies primitivity, an immediate consequence of Theorem~\ref{dppg-dcf} is that the same classification holds for differential-algebraic groups:

\begin{theorem}
\label{knop}
Suppose $(G,S)$ is a finite rank faithful $2$-transitive group action definable in $\dcf_0$ with $G$ connected.
Then $(G,S)$ is definably isomorphic to
\begin{itemize}
    \item[(i)]  $\mathrm{PSL}_{n+1}(\mathcal{C})$ acting on $\mathbb{P}^n(\mathcal{C})$, or
    \item[(ii)] $\mathbb G_a^n(\calC)\rtimes L$ acting on $\calC^n$ where $L=\mathrm{GL}_n(\calC)$ or $L=\mathrm{SL}_n(\calC)$, or
    \item[(iii)] $\mathbb G_a^{2n}(\calC)\rtimes L$ acting on $\calC^{2n}$ where $L=\mathrm{Sp}_{2n}(\calC)$ or $L=\mathrm{Sp}_{2n}(\calC)\cdot\mathbb G_m(\calC)$.
\end{itemize}
\end{theorem}

\begin{proof}
Note that a $2$-transitive group action is primitive (and hence definably primitive).
Indeed, if $E$ is an equivalence relation on $S$, and $x,y,z\in S$ are distinct elements such that $xEy$ but $\neg(xEz)$, then, as there is an element of $G$ taking the pair $(x,y)$ to $(x,z)$, it must be that $E$ is not $G$-invariant.
Hence, there are no proper nontrivial $G$-invariant equivalence relations.

Hence, by Theorem~\ref{dppg-dcf}, we have that $G$, and hence $(G,S)$, is definably isomorphic to a connected $2$-transitive algebraic group action in the constants, and these are classified by Knop~\cite{knop1983mehrfach} to be precisely of types~(i), (ii), or~(iii).
\end{proof}

There are many more {\em generically} $2$-transitive actions, even among algebraic group actions, than those appearing in Theorem~\ref{knop}.
However, a conjecture of Borovik and Cherlin from~\cite{BC2008} suggests that a high degree of generic transitivity is rare:

\begin{conjecture}[Borovik-Cherlin]
Suppose $(G,S)$ is a connected definable permutation group of finite Morley rank with $\mr(S)=n>0$.
If the action is generically $(n+2)$-transitive then $(G,S)$ is isomorphic to the natural action of $\mathrm{PSL}_{n+1}(F)$ on~$\mathbb{P}^n(F)$, for some algebraically closed field~$F$.
\end{conjecture}

In $\acf_0$, so for algebraic group actions, this conjecture was verified in~\cite{nmdeg}.
Using our Theorem~\ref{dppg-dcf} above, we will verify the conjecture in $\dcf_0$.
Note, however, that something more has to be done as generic $\mu$-transitivity (unlike outright $2$-transitivity) does not imply primitivity, and so Theorem~\ref{dppg-dcf} does not automatically apply.
Our proof will follow the suggestions in~\cite{BC2008} about how to reduce to the definably primitive case.

We work in a saturated model $(\calU,\delta)\models\dcf_0$ with field of constants~$\calC$.

\begin{theorem}[Borovik-Cherlin in $\dcf_0$]
\label{bc}
Suppose $(G,S)$ is a definable permutation group in $\dcf_0$ with $G$ connected and of finite rank, and $n=\mr(S)>0$.
If the action is generically $(n+2)$-transitive then $(G,S)$ is definably isomorphic to the natural action of $\mathrm{PSL}_{n+1}(\mathcal{C})$ on~$\mathbb{P}^n(\mathcal{C})$.
\end{theorem}

\begin{proof}
In~\cite[page~35]{BC2008} an argument is sketched for how to reduce the Borovik-Cherlin conjecture to definably primitive group actions.
In what follows we fill in some details, and implement some simplifications in the context of $\dcf_0$.

Let $H\leq G$ be a proper definable subgroup such that $S=G/H$.
First of all, notice that $H$ must be infinite.
Indeed, generic $(n+2)$-transitivity implies that $\mr(G)\geq n(n+2)$, and hence $\mr(G)>n$.

We proceed by induction on $n$.
What the induction hypothesis gives us is that if~$H'$ is any proper definable subgroup of $G$ containing $H$ then we must have $H'/H$ finite.
Indeed, consider the action of $G$ on $G/H'$ by left multiplication.
Note that it is also generically $(n+2)$-transitive as $H\leq H'$.
If we let $K\leq H'$ be the kernel of this action, then $(G/K,G/H')$ is faithful and generically $(n+2)$-transitive.
If $H'/H$ were infinite then $\mr(G/H')=:e<n$, and we can apply our induction hypothesis.
(In the case that $n=1$ we know by connectedness of $G$ that this does not happen, which deals with the base case.)
That is, we know that $(G/K,G/H')$ is definably isomorphic to the action of $\PGL_{e+1}(\calC)$ on $\P_e(\calC)$.
In particular, the action of $G$ on $G/H'$ is not generically $(e+3)$-transitive, contradicting the fact that the action is generically $(n+2)$-transitive and $e<n$.
Hence, we must have that $H'/H$ is finite.

Let $L$ be the normaliser of $H^\circ$ in $G$, where $H^\circ$ denotes the connected component of $H$.
We claim that $L$ is a finite extension of $H$.
Note that $$L:=\{g\in G:gH^\circ g^{-1}=H^\circ\}.$$
As such it is clearly a definable subgroup of $G$ containing $H$.
If $L=G$ then $H^\circ\lhd G$, and hence $H^\circ$ stabilises every point of $G/H$.
As the action of $G$ on $G/H$ is assumed to be faithful, this would imply that $H^\circ=(1)$, contradicting that $H$ is infinite.
So $L\neq G$.
Hence, as explained above, the induction hypothesis forces $L/H$ finite.

Note that $L$ contains every proper definable subgroup of $G$ that contains $H$.
Indeed, suppose $H'$ is such. Then, as $H'/H$ must be finite, we have that $H^\circ=(H')^\circ$.
It follows that $H^\circ\lhd H'$, forcing $H'\leq L$, as desired.

In particular, $L$ is a maximal proper definable subgroup of $G$.
So the action of $G$ on $G/L$ is definably primitive.
But it may no longer be faithful.
Let us show that the kernel of the action of $G$ on $G/L$, say~$K$, is finite. First note that $K = \bigcap\limits_{g \in G} L^g$. Similarly, because $G$ acts faithfully on $G/H$, we have $\bigcap\limits_{ g \in G} H^g = \{ 1 \}$.
By the descending chain condition on definable subgroups, there must be $g_1 , \cdots , g_n \in G$ such that $K = \bigcap\limits_{i = 1}^n L^{g_i}$ and $\{ 1 \} = \bigcap\limits_{i=1}^n H^{g_i}$. Note that for any subgroups $H_1 < L_1 < G$ and $H_2 < L_2 < G$, we have an injective map $(L_1 \cap L_2)/(H_1 \cap H_2)  \rightarrow L_1/H_1 \times L_2/H_2$. Applying this fact repeatedly, we see: 
\begin{align*}
    \left| K \right| & = \left[ K : \{ 1 \} \right] \\
    & = \left[  \bigcap\limits_{i = 1}^n L^{g_i}  :  \bigcap\limits_{i=1}^n H^{g_i}  \right] \\
    & \leq \prod\limits_{i = 1}^n \left[ L^{g_i} : H^{g_i} \right] \\
    & = \left[ L : H\right]^n
\end{align*}
\noindent and since $\left[ L : H \right]$ is finite, so is $K$.

We have that $(G/K,G/L)$ is faithful, transitive, and definably primitive.
Hence, by Theorem~\ref{dppg-dcf}, $(G/K,G/L)$ is definably isomorphic to an algebraic group action in the constants.
Note that $(G/K,G/L)$ is still generically $(n+2)$-transitive.
Now, the Borovik-Cherlin conjecture for $\acf_0$, as established in~\cite[Section~6]{nmdeg}, implies that $G/K$ is definably isomorphic to $\PGL_{n+1}(\calC)$.
That is, the quotient of $G$ by a normal finite subgroup is definably isomorphic to the $\calC$-points of an algebraic group over~$\calC$.
By~\cite[Corollary 3.10]{pillay-tehran}, this forces $G$ itself to be definably isomorphic to the $\calC$-points of an algebraic group over~$\calC$.
But then, $(G,S)$ is definably isomorphic to an algebraic group action in the constants, and the Borovik-Cherlin conjecture for $\acf_0$ implies it is definably isomorphic to $\mathrm{PSL}_{n+1}(\mathcal{C})$ on~$\mathbb{P}^n(\mathcal{C})$.
\end{proof}

The very same proof gives a differential-algebraic-geometric variant. In this variant we replace Morley-rank-based generic transitivity with the following natural notion of generic transitivity coming from the Kolchin topology:
 we say that a differential-algebraic group action, $(G,S)$, is {\em Kolchin-generically $\mu$-transitive} if the diagonal action of $G$ on $S^\mu$ admits a Kolchin dense orbit.
At the same time, we replace Morley rank itself with differential-algebraic-geometric dimension, which naturally generalises algebro-geometric dimension except that it is not always finite:
if $X$ is a differential-algebraic variety over a differential subfield $k$, then we say that $X$ is {\em finite dimensional} if, for each $a\in X$, the differential field $k\langle a\rangle$ generated by $a$ over $k$ is of finite transcendence degree over $k$, and in that case we call  the supremum of these transcendence degrees the {\em dimension} of $X$.

\begin{theorem}[Finite dimensional geometric BC for $\mathrm{DCF}_0$]
\label{fdgbc}
Suppose $(G,S)$ is a differential-algebraic permutation group with $G$ connected and finite dimensional, and $d=\dim S>0$.
If the action is Kolchin-generically $(d+2)$-transitive then $(G,S)$ is isomorphic to the natural action of $\mathrm{PSL}_{d+1}(\mathcal{C})$ on~$\mathbb{P}^d(\mathcal{C})$.
\end{theorem}

The precise relationship between Theorem~\ref{fdgbc} and Theorem~\ref{bc} is not clear.
First of all, finite-dimensionality is equivalent to finite Morley rank (in this single derivation case), and dimension is an upper bound for Morley rank.
However, dimension may be strictly larger than Morley rank.
Moreover, it is unlikely that the existence of a Kolchin dense orbit will coincide with that of a generic orbit in the sense of Morley rank.
For example, in his thesis~\cite[$\S$9]{freitag-thesis}, the first named author showed that the irreducible differential-algebraic variety $V$ defined by the ordinary order~3 algebraic differential equation $xx'''-x''=0$ has a proper Kolchin closed subset $W$ (defined by $x''=0$) whose complement has Morley rank strictly less than that of $V$.
So $W$ is generic in $V$ in the sense of Morley rank but is not Kolchin dense, while $V\setminus W$ is Kolchin dense in $V$ but not generic in the sense of Morley rank.
We do not, however, know of an example that arises as the orbit of a definable group action.

In any case, Theorem~\ref{fdgbc} follows from Theorem~\ref{dppg-dcf} in exactly the same way as Theorem~\ref{bc} did.
We leave the verification of details to the reader.

But our real purpose in raising this geometric variant is that it suggests (to us) a more general conjecture about $\dcf_{0,m}$, the theory of characteristic zero differentially closed fields in $m$ commuting derivations, $\Delta=\{\delta_1,\dots,\delta_m\}$.
The idea is that we can drop finite-dimensionality by replacing dimension with the Kolchin polynomial, a numerical polynomial associated to a differential-algebraic variety measuring the growth in transcendence degree as you take higher order derivatives of a solution.
See the exposition in~\cite{cassidy2011jordan} for details.
The {\em $\Delta$-type} is the degree of the Kolchin polynomial (which is at most $m$) and the {\em typical $\Delta$-dimension} is its leading coefficient.
The finite dimensional case corresponds to when the $\Delta$-type is zero, and in that case the typical $\Delta$-dimension is what we have been calling dimension.
Instead of connectedness we have to consider the following strengthening introduced in \cite{cassidy2011jordan}:
A differential-algebraic group $G$ is {\em strongly connected} if there is no proper differential-algebraic subgroup $H$ such that $\deltatype(G/H) < \deltatype(G)$.
Strongly connected groups are connected, and the converse also holds of finite dimensional groups.
Here is what Theorem~\ref{fdgbc} leads us to expect:

\begin{conjecture}[Geometric BC for $\mathrm{DCF}_{0,m}$]
\label{gbc}
Suppose $G$ is a strongly connected differential-algebraic group of $\Delta$-type $\ell$ acting differential-algebraically, faithfully and transitively on a differential-algebraic variety $S$ of typical $\Delta$-dimension $d>0$.
If $G$ acts Kolchin-generically $(d+2)$-transitively on $S$ then 
$(G,S)$ is isomorphic to the natural action of $\mathrm{PSL}_{d+1}(\calF)$ on~$\mathbb{P}^d(\calF)$ where $\calF$ is the constant field of some $m-\ell$ linearly independent derivations in $\Span_{\calC}(\Delta)$.
\end{conjecture}

Theorem~\ref{fdgbc} is the $\deltatype(G)=0$ case of Conjecture~\ref{gbc}.
It is the opposite extreme, when $\deltatype(G)=m$, that poses the greatest difficulty.
If this case of maximal $\Delta$-type were settled, then an induction on the number of derivations, using the work of Le\'on S\'anchez~\cite{ls15} on {\em relative} $D$-groups and $D$-varieties, together with Buium's theorem on the isotriviality of $D$-variety structures on projective varieties~\cite{buium2007differential}, should imply Conjecture~\ref{gbc}.
However, the maximal $\Delta$-type case, even when $m=1$, remains open.

\bigskip
\section{Bounding nonorthogonality}
\label{nonorth}

\noindent
The permutation groups in $\dcf_0$ that we are primarily interested in are those that arise as binding groups.
The starting point of~\cite{nmdeg} was the observation that bounding the generic transitivity degree of binding group actions leads to a bound on the witness to nonorthogonality.

Recall that a type $p\in S(k)$ is {\em weakly orthogonal}\footnote{This is also called ``almost orthogonal" in the literature.} to a $k$-definable set $X$ if every realisation of $p$ is independent of every finite tuple from $X$ over $k$.
It is outright {\em orthogonal} if this holds of all nonforking extension of $p$ to more parameters.
It is a general fact that if $p$ is nonorthogonal to $X$ then some  finite Morley power of $p$ is not weakly orthogonal to $X$.
A fundamental theorem of Hrushovski's implies that
{\em if $p$ is minimal and nonorthogonal to~$\calC$ then $p^{(4)}$ is already not weakly orthogonal to~$\calC$}.
Indeed this is a consequence of a very special case of~\cite[Theorem~1]{hrushovski1989almost}, which is about arbitrary stable theories and (possibly infinite rank) regular types.
One cannot expect in $\dcf_0$ an absolute bound that is independent of rank: the construction in~\cite[Section~4.2]{nmdeg} exhibits, for each $n\geq 2$, a Morley rank $n$ type $p$ that is nonorthogonal to $\calC$ (indeed, is $\calC$-internal), and such that $p^{(n+2)}$ is weakly $\calC$-orthogonal.\footnote{However, if we assume in addition that $\delta$ is trivial on $k$, then it is shown in~\cite{abred} that $p$ being nonorthogonal to $\calC$ implies $p^{(2)}$ is already not weakly $\calC$-orthogonal.}
However, this is as bad as it gets; applying the truth of the Borovik-Cherlin conjecture to binding groups, the first and third authors showed in~\cite{nmdeg} that $n+3$ is a bound for witnessing nonorthogonality to the constants.
Using Borovik-Cherlin in $\dcf_0$, namely Theorem~\ref{bc} above, instead, the proof extends to arbitrary definable sets:

\begin{corollary}
\label{cor:int}
Suppose $p\in S(k)$ is a type of Morley rank $n$ over an algebraically closed differential field $k$, and $X$ is any definable set over $k$.
If $p$ is nonorthogonal to $X$ then the Morley power $p^{(n+3)}$ is not weakly $X$-orthogonal.

When $X=\calC$ we can replace Morley rank with $U$-rank.
\end{corollary}

\begin{proof}
From nonorthogonality to $X$ we get a definable function $p\to q$ such that $q$ is nonalgebraic and $X$-internal (see~\cite[7.4.6]{GST}).
So we may as well assume that $p$ is $X$-internal.
Let $G$ be the binding group of $p$ relative to $X$.
As we may assume that $p$ itself is weakly $X$-orthogonal, $G$ acts transitively on $S:=p(\calU)$.
In particular, $p$ is isolated and $S$ is a definable set of Morley rank $n$.
(Note that, as $p$ is $X$-internal, if $X=\calC$ then $n=U(p)$.)
Moreover, as $k$ is algebraically closed, $G$ is connected.
If $p^{(n+3)}$ is weakly $X$-orthogonal then $G$ acts transitively on $p^{(n+3)}(\calU)$, and so $(G,S)$ is generically $(n+3)$-transitive.
But then it is also generically $(n+2)$-transitive, and so, by our Theorem~\ref{bc}, $(G,S)$ is $(\PGL_{n+1}(\calC),\mathbb P^n(\calC))$.
But this is a contradiction as the latter is {\em not} generically $(n+3)$-transitive.
\end{proof}

The above Corollary answers, for the case of finite rank types in $\dcf_0$, a question Hrushovski raises in~\cite[Section~2]{hrushovski1989almost} as a possible refinements of his main theorem; it is the second of the two ``natural generalizations" that he mentions there.

In fact, with a little more work, we can obtain a more refined version of Corollary~\ref{cor:int}.
So far we have been talking about definable interaction between single realisations of $p$ and finite tuples from $X$.
But one can ask the finer, more symmetric question, about nonorthogonality between two types: $p$ and $q$ are {\em weakly orthogonal}, denoted $p\perp^w q$, if any realisation of $p$ is independent of any realisation of $q$, and they are {\em orthogonal}, denoted $p\perp q$, if this continues to hold after taking nonforking extensions.
If $p$ and $q$ are nonorthogonal then some finite Morley powers of $p$ and $q$ will be weakly nonorthogonal.
Can we bound these powers?
This question, for arbitrary stable theories but working with regular types, was also raised by Hrushovski in~\cite{hrushovski1989almost}; it is the other of the two natural generalizations he asks about.
We show:

\begin{theorem}
\label{orthbound}
Suppose $k$ is an algebraically closed differential field and $p,q \in S(k)$ are of $U$-rank $m$ and $n$ respectively.
If $p \not\perp q$ then $p^{(m+3)} \not\perp^w q^{(n+3)}$.
\end{theorem}

Before proving the theorem we recall two well known facts from geometric stability theory, including proofs for the sake of completeness.
First, all instances of nonorthogonality between finite rank types are witnessed by minimal types:

\begin{lemma}\label{lem:minmid}
Suppose $p$ and $q$ are complete stationary types of finite rank.
Then $p\not\perp q$ if and only if there is a minimal type $r$ such that $p\not\perp r$ and $r\not\perp q$.
\end{lemma}

\begin{proof}
For the right-to-left direction (which does not use finite rank), suppose $p$ and $q$ are each nonorthogonal to some minimal $r$.
Taking nonforking extensions we may assume that in fact all three types are over a common parameter set $A$, and that $p$ and $q$ are in fact each non-weakly-orthogonal to $r$.
Then there are $a\models p, b\models q, c\models r$ such that $a\nind_A c$ and $c\nind_A b$.
Since $r$ is minimal, we have $c\in\big(\acl(aA)\cap\acl(bA)\big)\setminus\acl(A)$ witnessing that $a\nind_Ab$.

For the converse, we use \cite[Corollary~1.4.5.7]{GST}: for any finite rank $p \in S(A)$, there exists a model $A \subset M \models T$ and minimal types $p_1 \cdots , p_n \in S(M)$ such that $p_M$ is {\em domination equivalent} to $p_1 \otimes \cdots \otimes p_n$.
    That is, a tuple forks over $M$ with a realisation of $p_M$ if and only if it forks over $M$ with a realisation of $p_1 \otimes \cdots \otimes p_n$.
    Taking nonforking extensions to a larger model if necessary, we may assume that $M$ contains the domain of $q$, and that, since $p\not\perp q$, a realisation of $q_M$ forks with a realisation of $p_M$ over $M$.
    Forking calculus now shows that $p_i \not\perp q$, for some $i\leq n$.
    As $p \not\perp p_i$, we can take $r=p_i$.
\end{proof}

Next, recall that for a stationary type $p\in S(A)$ and a set of complete types $\calQ$ over arbitrary parameter sets, we say that $p$ is {\em $\calQ$-internal} if there exist $\ell\geq 1$,  $\tp(b_i/B_i)\in\calQ$ for
$i=1,\dots,\ell$, $C\supseteq A\cup \bigcup_{i=1}^\ell B_i$, and $a\models p$ with $a\ind_AC$, such that $a\in\acl(Cb_1,\dots,b_\ell)$.
In particular, for $q\in S(B)$, we have that $p$ is {\em $q$-internal} if it is $\{q\}$-internal.
Replacing $\dcl$ with $\acl$ yields the notion of {\em almost internality}.

Our second lemma says that almost internality to a minimal type coincides with almost internality to the set of conjugates of that minimal type:

\begin{lemma}
\label{lem: int-conj}
Suppose $p\in S(A)$ and $q\in S(B)$ are stationary with $q$ minimal and where $B\supseteq A$.
Let $\calQ$ be the set of $A$-conjugates of $q$.
If $p$ is almost $\calQ$-internal then $p$ is almost $q$-internal.
\end{lemma}

\begin{proof}
(Thanks to Anand Pillay for pointing out this argument.)

If $p=\tp(a/A)$ is almost $\calQ$-internal then there are $A$-conjugates of $q$, say $q_i=\tp(b_i/B_i)$, $i=1,\dots,\ell$, and $C\supseteq \bigcup_{i=1}^\ell B_i$ with $a\ind_AC$, such that $a\in\acl(Cb_1,\dots,b_\ell)$.
Choose these so that $\ell$ is minimal.

If $b_1\in\acl(C)$ then we can just drop it and contradict minimality of $\ell$.
So we may assume $b_1\notin\acl(C)$.
If $b_1\notin\acl(Ca)$ then $a\ind_Cb_1$ as $q$ is of $U$-rank~$1$, and replacing~$C$ with $Cb_1$ we contradict the minimality of $\ell$.
Hence $b_1\in\acl(Ca)\setminus\acl(C)$.
It follows that $q_1=\tp(b_1/B_1)$ is nonorthogonal to $p=\tp(a/A)$.
This means (by definition) that $q_1$ is nonorthogonal to the parameter set $A$.
By \cite[Lemma 1.4.3.3]{GST}, we get that if $B'$ is an $A$-conjugate of $B$ that is independent from $B_1$ over $A$, and $q'$ is the corresponding $A$-conjugate of $q$, then $q_1$ is nonorthogonal to $q'$.
Of course this holds not just for $q_1$ but for $q_2,\dots,q_\ell$ as well.
So if we choose $B'$ to be an $A$-conjugate of $B$ that is independent of $(B_1,\dots,B_\ell)$ over $A$, and let $q'$ be the corresponding conjugate of $q$, then we have that each $q_i$ is nonorthogonal to~$q'$.
As these are minimal types we get that $q_1,\dots,q_\ell$ are pairwise nonorthogonal.

So, taking a larger $C$ and further nonforking extension if necessary, we can assume that each $b_i$ is interalgebraic with some $c_i$ over $C$, where all the $c_1,\dots,c_\ell$ realise $q_1$.
Hence $a\in\acl(Cc_1,\dots,c_\ell)$ witnesses almost internality of $p$ with $q_1$.
As $q_1$ is an $A$-conjugate of $q$, we have that $p$ is almost $q$-internal.
\end{proof}

We are now ready to prove the theorem.

\begin{proof}[Proof of Theorem~\ref{orthbound}]
    By Lemma \ref{lem:minmid}, there exists a minimal type $t$, potentially over additional parameters, such that $p \not\perp t$ and $q \not \perp t$.
    We deal separately with the cases when $t$ is or is not locally modular.

    Suppose first that $t$ is non-locally-modular.
    Then both $p$ and $q$ are nonorthogonal to the constants.
    By Corollary~\ref{cor:int} we have that $p^{(m+3)} \not\perp^w \mathcal{C}$.
    It follows that there is a realisation $a\models p^{(m+3)}$ and a $k$-definable function, $f$, such that $f(a)\in\calC$ is generic over $k$
    (see, for example, Lemma~2.1 of~\cite{nmdeg}.)
    Similarly, we have $b\models q^{(n+3)}$ and a $k$-definable function $g$, such that $g(a)\in\calC$ is generic over $k$.
    Since $f(a)$ and $g(b)$ realise the same type over $k$, by automorphisms, we can choose $b$ such that $f(a)=g(b)$.
    It follows that $a\nind_kb$, witnessing that $p^{(m+3)} \not\perp^w q^{(n+3)}$.

    Now suppose that $t$ is locally modular.
    We show that in this case $p^{(2)}\not\perp^wq^{(2)}$.
    
    Let $\mathcal{T}$ be the set of $k$-conjuguates of $t$.
    Since $p$ and $q$ are nonorthogonal to $t$, we obtain, by~\cite[7.4.6]{GST}, $k$-definable functions $p\to r$ and $q\to s$ such that $r$ and $s$ are nonalgebraic and $\mathcal{T}$-internal.
    Hence $r$ and $s$ are almost $t$-internal by Lemma~\ref{lem: int-conj}.
    In particular, $r$ and $s$ are $1$-based. At this point, given $r$ and $s$ one-based and nonorthogonal, it is well known how to deduce that $r^{(2)}$ and $s^{(2)}$ are non-weakly orthogonal. Nevertheless, we give some details.
    
    We claim that there is a minimal type over $k$ that is algebraic over $r$.
    That is, writing $r=\tp(c/k)$, there is $\widetilde c\in\acl(kc)$ such that $\widetilde r:=\tp(\widetilde c/k)$ is minimal.
    If $U(r)=1$ then we can simply take $\widetilde c=c$.
    Otherwise, there is $B \supset k$ with $U(c/B) = U(c/A) - 1>0$.
    Let $\widetilde c = \mathrm{Cb}(c/B)$.
    By $1$-basedness,  $\widetilde c \in \acl(kc)$.
    We have $U(c/k\widetilde c) = U(c/k) - 1$, from which it follows that $U(\widetilde c/k) = 1$, as desired.
    Similarly, we have $s=\tp(d/k)$ and $\widetilde d\in\acl(kd)$ such that $\widetilde s:=\tp(\widetilde d/k)$ is minimal.
    
    Note that $\widetilde r$ and $\widetilde s$ are locally modular nonorthogonal minimal types.
    Now, for modular minimal types orthogonality and weak orthogonality coincide, see~\cite[2.5.5]{GST}.
    It follows that $\widetilde r^{(2)} \not\perp^w \widetilde s^{(2)}$, and hence $r^{(2)} \not\perp^w s^{(2)}$.
    As $r$ and $s$ are $k$-definable images of $p$ and $q$, we have that $p^{(2)} \not\perp^w q^{(2)}$, as desired.
\end{proof}

\begin{remark}
 While the proof of Theorem~\ref{orthbound} appealed to Corollary~\ref{cor:int}, it only used the $X=\calC$ case of that corollary, and hence only makes use of the truth of Borovik-Cherlin in $\acf_0$.
\end{remark}

\bigskip
\section{Other theories}
\label{ot}

\noindent
So far we have worked in the theory $\dcf_0$.
Here we verify that the results hold also, for the finite rank setting, in $\dcf_{0,m}$ and the theory of compact complex manifolds ($\ccm$).
In fact we extract from the above proofs the required abstract properties of the theory.

\begin{theorem}
\label{dppg-T}
Suppose $T$ is a complete stable theory admitting elimination of imaginaries, and with a pure $0$-definable algebraically closed field $\calC$.
Assume the following hold in $T$:
\begin{itemize}
    \item[(1)] Every simple group of finite rank definable in $T$ is definably isomorphic to the $\calC$-points of a simple linear algebraic group over~$\calC$.

    \item[(2)] Up to definable isomorphism, the only commutative finite rank group definable in $T$, having no proper nontrivial definable subgroups, is $\mathbb G_a(\calC)$.

    \item[(3)] There are no infinite elementary abelian groups of finite rank definable in~$T$.
\end{itemize}
Then Theorems~\ref{A} and~\ref{B} of the Introduction hold of $T$, with $\calC$ playing the role of the constants.

Suppose, moreover, that the following holds in $T$:
\begin{itemize}
    \item [(4)]
    Every connected commutative finite rank group definable in $T$ is divisible.
\end{itemize}
Then the Borovik-Cherlin Conjecture, namely Theorem~\ref{C}, holds in $T$, again with $\calC$ playing the role of the constants.

In particular, Theorems~\ref{A}, \ref{B}, \ref{C} of the Introduction all hold of $\dcf_{0,m}$ with $\calC$ the field of total constants, and of $\ccm$ with $\calC$ the elementary extension of the complex field living on the projective line.\footnote{Note that the Borovik-Cherlin Conjecture in $\ccm$ was already established in~\cite{nmdeg}.}
\end{theorem}

\begin{proof}
We leave it to the reader to verify that conditions~(1) through~(3) are exactly what is needed for the proof of Theorem~\ref{A} (namely~\ref{dppg-dcf}) to go through.
Theorem~\ref{B} (namely~\ref{bd-dcf}) for $T$ follows formally from Theorem~\ref{A} for $T$ together with Theorem~\ref{bd-acf} (for $\acf_0$).

Maybe it is worth pointing out that condition~(1) already implies that, up to definable isomorphism, $\calC$ is the only infinite field of finite rank definable in $T$.
See, for example, Pillay's proof of Corollary~1.6 in~\cite[Chapter~III]{marker2017model}.
Namely, if $F$ is an infinite definable field of finite rank in $T$, then $\PGL_2(F)$ is a simple group of finite rank definable in $T$, and hence $\PGL_2(F)$ is definably isomorphic to a group definable in the pure field $\calC$.
But the field $F$ itself is interpretable in the group structure on $\PGL_2(F)$, so that $F$ is then definably isomorphic to a field definable in $(\calC,+,\times)$, and hence is definably isomorphic to $\calC$.

It is shown in~\cite[Corollary~3.10]{pillay1997some} that condition~(4) implies that if $G$ is a connected definable group with the property that the quotient by some normal finite subgroup is definably isomorphic to the $\calC$-points of an algebraic group over~$\calC$ then this is already true of $G$.
With this additional fact, our proof of Theorem~\ref{bc} goes through to establish the Borovik-Cherlin conjecture in $T$.

Finally, let us point out that these conditions do hold in $\dcf_{0,m}$ and $\ccm$.

In $\dcf_{0,m}$, as in the case of a single derivation, it remains true that every definable group of finite rank embeds in an algebraic group.
Conditions~(3) and~(4) follow, as every commutative algebraic group in characteristic zero has finite $n$-torsion for all $n$.
From this embedding into an algebraic group it also follows that if $G$ is a simple definable group of finite rank, then it embeds definably in $\GL_n$ (see, for example, the proof of Corollary~1.2 in~\cite[Chapter~III]{marker2017model}), and hence Cassidy's theorem (see~\cite[Theorem~3.7]{cassidy2011jordan} for a version that holds also of infinite rank) implies that $G$ is definably isomorphic to the $\calC$-points of a simple linear algebraic group.
So condition~(1) holds.
For condition~(2), note that the proof of Lemma~\ref{nosbgp} above goes through in the case of $\dcf_{0,m}$ (in particular, \cite{hrushovskiminimal} works in the context of several commuting derivations).

Regarding $\ccm$, there is a very good understanding of the definable groups coming out of~\cite{fujiki-mergroups, pillay-scanlon-mergroups, ams-mergroups, scanlon-mergroups}.
In particular, there is a Chevalley-type structure theorem for definable groups whereby every definable group is the extension of a ``nonstandard complex torus" by a linear algebraic group.
Here the notion of a nonstandard complex torus is somewhat subtle, but at the very least these are commutative groups with nontrivial but finite $n$-torsion for all $n$.
Conditions~(1) through~(4) all follow from this.
\end{proof}

\begin{theorem}
\label{orthbound-T}
Suppose $T$ is a complete totally transcendental theory admitting elimination of imaginaries, and with a pure $0$-definable algebraically closed field~$\calC$, such that every non-locally-modular minimal type is nonorthogonal to $\calC$.
Then Theorem~\ref{D} of the Introduction holds in $T$.

In particular, this holds of $\dcf_{0,m}$ and $\ccm$.
\end{theorem}

\begin{proof}
Again, we leave to the reader the verification that these were the only properties of $\dcf_0$ used in the proof of Theorem~\ref{D} (namely, \ref{orthbound}).
That they hold of $\dcf_{0,m}$ and $\ccm$ is well known.
\end{proof}

\bigskip

\bibliography{biblio}
\bibliographystyle{plain}
\end{document}